\documentclass[12pt, reqno]{amsart}
\usepackage[margin=1in]{geometry}

\usepackage{mathtools}
\usepackage{amssymb, amsthm, thmtools, mleftright}
\usepackage[backend=biber, style=numeric, sorting=nyt, url=true, doi=true]{biblatex}
\usepackage[stretch=10, shrink=10]{microtype}
\usepackage{booktabs, multirow}
\usepackage{xurl}
\usepackage{hyperref}
\DeclareFieldFormat{doi}{%
  \newline
  \mkbibacro{DOI}\addcolon\space
    \ifhyperref
      {\href{http://dx.doi.org/#1}{\nolinkurl{#1}}}
      {\nolinkurl{#1}}}

\DeclareMathOperator{\res}{Res}
\NewDocumentCommand{\abs}{m}{\mleft\lvert #1 \mright\rvert}
\NewDocumentCommand{\norm}{m}{\mleft\| #1 \mright\|}
\NewDocumentCommand{\df}{}{\mathop{}\!\mathrm{d}}

\NewDocumentCommand{\numberthis}{}{\addtocounter{equation}{1}\tag{\theequation}}

\AtBeginDocument{\RenewDocumentCommand{\Re}{}{\operatorname{Re}}}
\AtBeginDocument{\RenewDocumentCommand{\Im}{}{\operatorname{Im}}}
\AtBeginDocument{}

\numberwithin{equation}{section}
\declaretheorem[numberwithin=section]{theorem}
\declaretheorem[numberlike=theorem]{lemma, corollary, definition, proposition}

\begin{document}

\title{Euler Product Sieve}

\author{Di Liu}
\address{Department of Mathematics, University of Illinois at Urbana--Champaign,
1409 West Green Street, Urbana, IL 61801, USA}
\email{dil4@illinois.edu}

\author{ {Yuri} Matiyasevich}
\address{St.~{Petersburg} Department of Steklov Mathematical Institute
of Russian Academy of Sciences, nab. Fontanki 27, St. Petersburg, 191023 Russia}
\email{yumat@pdmi.ras.ru}

\author{Joseph Oesterl\'{e}}
\address{ {Institut de Math\'{e}matique de Jussieu--Paris Rive Gauche,
4 place Jussieu,
Boite Courrier 247,
75252 Paris Cedex 5, France}}
\email{ {joseph.oesterle@imj-prg.fr}}

\author{Alexandru Zaharescu}
\address{Department of Mathematics, University of Illinois at Urbana--Champaign,
1409 West Green Street, Urbana, IL 61801, USA \and
Simion Stoilow Institute of Mathematics of the Romanian Academy,
P.O. Box 1-764, RO-014700, Bucharest, Romania}
\email{zaharesc@illinois.edu}

\date{}

\begin{abstract}
   {
  We study a class of approximations to
  the Riemann zeta function introduced
  earlier by the second author on the basis of
  Euler product.
  This allows us to justify
  \emph{Euler Product Sieve}
  for generation of prime numbers.}

   {Also we show that \emph{Bounded
  Riemann Hypothesis} (stated in a paper by
  the fourth author) is equivalent
  to conjunction \emph{the Riemann Hypothesis + the simplicity of zeros}.}
\end{abstract}



\maketitle

\section{Introduction}

In~\cite{matiyasevichFewFactorsEuler2017}, the second author constructed
a novel family of approximations to the Riemann zeta and xi functions
denoted by \( \zeta_{\mathfrak{B}}^{\approx}\mleft( s \mright) \)~
and~\( \xi_{\mathfrak{B}}^{\approx}\mleft( s \mright) \),
the precise definitions of which are given in~\eqref{approx} below.
Many interesting questions were raised in~\cite{matiyasevichFewFactorsEuler2017}.
Some of them were answered by Nastasescu and the fourth author
in~\cite{nastasescuClassApproximationsRiemann2022}.
A striking question not answered in~\cite{nastasescuClassApproximationsRiemann2022} concerns
an original sieve for prime numbers,
called the \emph{Euler product sieve}.
The sieve is defined as follows,
for more details, see~\cite[Section~4]{matiyasevichFewFactorsEuler2017}.

\begin{definition}
  \label{def:sieve}
  Let the function~\( \xi_{\mathfrak{B}}^{\approx}\mleft( s \mright) \)
  be defined by~\eqref{approx}. For a given complex number~\( s_{0} \),
  Euler product sieve at~\( s_{0} \) is the following procedure to locate primes:
  \begin{enumerate}
    \item Announce a few initial prime numbers~\( 2, \ldots, p_{k_0} \).
    \item Suppose the first \( k \)~primes \( 2, \ldots, p_{k} \) have been found,
      calculate the difference
      \[
        \delta_{\mleft\{ 2, \ldots, p_{k} \mright\}} \coloneqq
        \xi_{\mleft\{ 2, \ldots, p_{k} \mright\}}^{\approx}\mleft( s_0 \mright) - \xi\mleft( s_0 \mright).
      \]
    \item The next prime is the least natural number~\( q \) which
    \begin{itemize}
      \item is greater than~\( p_{k} \), and different from the powers of~\( 2, \ldots, p_{k} \), and
      \item gives better approximation than the one using \( \mleft\{ 2, \ldots, p_{k} \mright\} \) only, i.e.,
      \[
        \abs{ \xi_{\mleft\{ 2, \ldots, p_{k}, q \mright\}}^{\approx}
        \mleft( s_0 \mright) - \xi\mleft( s_0 \mright) }
        < \abs{ \delta_{\mleft\{ 2, \ldots, p_{k} \mright\}} }.
      \]
    \end{itemize}
  \end{enumerate}
\end{definition}
\autoref{tbl:sieve} demonstrates how the sieve picks a few initial primes
for~\( s_{0} = 1/2 \).
In the present paper we prove that for~\( s_0 = 1/2 \) it is sufficient
to start by announcing just the first prime number.
\begin{table}[ht]
  \centering
  \caption{Euler product sieve}
  \begin{tabular}{ll}
    \toprule
    \( \mathfrak{B} \) & \( \xi_{\mathfrak{B}}^{\approx}\mleft( 1/2 \mright) - \xi\mleft( 1/2 \mright) \) \\
    \midrule
    \( \mathfrak{B}_{\le 7} \cup \mleft\{ 10 \mright\} \)
    & \( -3.73 \ldots \cdot 10^{-103} \) \\
    \( \mathfrak{B}_{\le 7} \cup \mleft\{ 11 \mright\} \)
    & \( 9.34 \ldots \cdot 10^{-195} \) \\[2mm]
    \( \mathfrak{B}_{\le 13} \cup \mleft\{ 14 \mright\} \)
    & \( -2.37 \ldots \cdot 10^{-217} \) \\
    \( \mathfrak{B}_{\le 13} \cup \mleft\{ 15 \mright\} \)
    & \( -5.22 \ldots \cdot 10^{-256} \) \\
    \( \mathfrak{B}_{\le 13} \cup \mleft\{ 17 \mright\} \)
    & \( 1.85 \ldots \cdot 10^{-438} \) \\[2mm]
    \( \mathfrak{B}_{\le 23} \cup \mleft\{ 24 \mright\} \)
    & \( 1.31 \ldots \cdot 10^{-705} \) \\
    \( \mathfrak{B}_{\le 23} \cup \mleft\{ 26 \mright\} \)
    & \( 1.41 \ldots \cdot 10^{-840} \) \\
    \( \mathfrak{B}_{\le 23} \cup \mleft\{ 28 \mright\} \)
    & \( 1.42 \ldots \cdot 10^{-986} \) \\
    \( \mathfrak{B}_{\le 23} \cup \mleft\{ 29 \mright\} \)
    & \( -3.34 \ldots \cdot 10^{-1226} \) \\
    \bottomrule
  \end{tabular}
  \label{tbl:sieve}
\end{table}

\begin{theorem}
  \label{thm:sieveAtOneHalf}
  The Euler product sieve at \( s = 1/2 \) always picks the primes
  after the first prime~\( 2 \) is announced.
\end{theorem}
In the case of general \( s_{0} \in \mathbb{C} \setminus \mleft\{ 0, 1 \mright\} \)
we might need to announce more than one initial primes.

\begin{theorem}
  \label{thm:sieve}
  For any \( s_{0} \in \mathbb{C} \setminus \mleft\{ 0, 1 \mright\} \),
  there exists a \( k_{0} \in \mathbb{N} \) depending only on~\( s_{0} \),
  such that if the first \( k_{0} \) primes are announced,
  then the Euler product sieve at \( s_{0} \) always picks the remaining primes.
\end{theorem}

A natural question following \autoref{thm:sieveAtOneHalf}
and \autoref{thm:sieve}, which we leave to interested readers, is to
characterize for each natural number~\( k \),
the set of points in the complex plane which require
the first \( k \) primes to be announced,
so that the sieve at those points always picks primes afterwards.
In order to state this question more clearly, we introduce the following family of sets.
\begin{definition}
  \label{def:dSet}
  For \( k \in \mathbb{N} \), let
  \begin{multline*}
    \mathcal{D}_{k} \coloneqq \mleft\{ s \in \mathbb{C} \setminus \mleft\{ 0, 1 \mright\}
    \colon \text{the Euler product sieve always picks primes,} \mright. \\
    \mleft. \text{once the first \( k \) primes are announced.} \mright\}.
  \end{multline*}
\end{definition}
By~\autoref{def:dSet}, we see the sets~\( \mleft( \mathcal{D}_{k} \mright) \)
form an ascending chain, i.e.,
\( \mathcal{D}_{1} \subseteq \mathcal{D}_{2} \subseteq \ldots \).
Moreover, \autoref{thm:sieveAtOneHalf} can be restated as \( 1/2 \in \mathcal{D}_{1} \),
while \autoref{thm:sieve} can be restated as
\( \cup_{k \in \mathbb{N}} \mathcal{D}_{k} = \mathbb{C} \setminus \mleft\{ 0, 1 \mright\} \).
We leave as an open question to interested readers,
the problem of characterzing the above chain of sets
\( \mathcal{D}_{1} \subseteq \mathcal{D}_{2} \subseteq \ldots\).

Apart from the Euler product sieve,
there are several other questions which can now be answered.
The first one is \cite[(5.1)]{matiyasevichFewFactorsEuler2017},
on the limiting behavior of those approximation functions.
\begin{theorem}
  \label{thm:limit}
  Fix an \( m \in \mathbb{N} \),
  let \( \zeta_{\mathfrak{B}}^{\approx}\mleft( s \mright) \) be as in~\eqref{approx},
  \( g\mleft( s \mright) \) be as in~\eqref{xi} and
  \begin{equation}
    \label{zeta_m}
    \zeta_{\mleft[ m \mright]}\mleft( s \mright)
    \coloneqq \lim_{\min \mleft\{ p_{1}, \ldots, p_{m} \mright\} \to \infty}
    \zeta_{\mleft\{ p_{1}, \ldots, p_{m} \mright\}}^{\approx}\mleft( s \mright).
  \end{equation}
  Then
  \[
    \zeta_{\mleft[ m \mright]}\mleft( s_{0} \mright)
    = \frac{1}{s_{0}^{m} \mleft( 1 - s_{0} \mright)^{m} g\mleft( s_{0} \mright)} \cdot
    \frac{1}{2 \pi i} \int_{\mleft( \abs{ \Re s_{0} } + 2 \mright)}
    s^{m} \mleft( 1 - s \mright)^{m}
    g\mleft( s \mright) a\mleft( s, s_{0} \mright) \df s,
  \]
  where \( a\mleft( s, s_{0} \mright)
  \coloneqq 1/\mleft( s_{0} - s \mright) + 1/\mleft( s_{0} - 1 + s \mright) \).
\end{theorem}
Another question is \cite[(A3)]{matiyasevichFewFactorsEuler2017},
concerning the sign alternation at composite numbers.
\begin{theorem}
  \label{thm:sign}
  Let \( m \in \mathbb{N} \)~and~\( \mathfrak{B} \) be the set of the first \( m \) primes.
  Let \( k \in \mathbb{N} \) be a composite number between \( p_{m} \)~and~\( p_{m + 1} \),
  then
  \[
    \Big( \xi_{\mathfrak{B}}^{\approx}\mleft( 1/2 \mright) - \xi\mleft( 1/2 \mright) \Big)
    \mleft( \xi_{\mathfrak{B} \cup \mleft\{ k \mright\}}^{\approx}\mleft( 1/2 \mright)
    - \xi\mleft( 1/2 \mright) \mright) > 0.
  \]
\end{theorem}

Finally we show the equivalence between the Bounded Riemann Hypothesis
\cite[Conjecture~4]{nastasescuClassApproximationsRiemann2022} and
Riemann Hypothesis plus the Simplicity Conjecture.

\begin{theorem}
  \label{thm:rhsc}
  The Bounded Riemann Hypothesis is equivalent to
  the Riemann Hypothesis plus the Simplicity Conjecture.
\end{theorem}

\section{Construction}

We first recall the construction of
several auxiliary functions
(for more details, see~\cite[Section~2]{matiyasevichFewFactorsEuler2017}).
Let \( \mathfrak{B} \) be a subset of the primes,
we define the finite Euler product~\( \zeta_{\mathfrak{B}}\mleft( s \mright) \) by
\begin{equation}
  \label{zeta}
  \zeta_{\mathfrak{B}}\mleft( s \mright)
  \coloneqq \prod_{p \in \mathfrak{B}} \frac{1}{1 - p^{-s}}.
\end{equation}
as well the completed version~\( \xi_{\mathfrak{B}}\mleft( s \mright) \) by
\begin{equation}
  \label{xi}
  \xi_{\mathfrak{B}}\mleft( s \mright)
  \coloneqq g\mleft( s \mright) \zeta_{\mathfrak{B}}\mleft( s \mright)
  = \tfrac{s\mleft( s - 1 \mright)}{2} \pi^{-s/2}
  \Gamma\mleft( s/2 \mright) \zeta_{\mathfrak{B}}\mleft( s \mright).
\end{equation}
Next we symmetrize \( \xi_{\mathfrak{B}}\mleft( s \mright) \)
by defining \( \xi_{\mathfrak{B}}^{=}\mleft( s \mright) \) as
\begin{equation}
  \label{eqxi}
  \xi_{\mathfrak{B}}^{=}\mleft( s \mright)
  \coloneqq \xi_{\mathfrak{B}}\mleft( s \mright) + \xi_{\mathfrak{B}}\mleft( 1 - s \mright),
\end{equation}
so now \( \xi_{\mathfrak{B}}^{=}\mleft( s \mright) \)
is equipped with a functional equation.
Now we consider some analytic properties of~\( \zeta_{\mathfrak{B}}\mleft( s \mright) \).
For each prime~\( p \in \mathfrak{B} \), \eqref{zeta} has poles located at
\begin{equation}
  \label{pPole}
  s = n \frac{2 \pi i}{\log p},
  \text{for~}
  n \in \mathbb{Z}.
\end{equation}
Those are called \emph{ghost poles of the zeta function}
by the second author, because they do not appear in the Riemann zeta function.
When \( n \ne 0 \), all of them are different simple poles with residue~\( 1/ \log p \).
However when \( n = 0 \), they all coincide and contribute to a pole of order~\( m \) at~\( s = 0 \)
where \( m = \abs{ \mathfrak{B} } \). One would want to remove those poles so the function has better
analytic properties, and this process is carried out below.

First we define a new function~\( \xi_{\mathfrak{B}}^{\coloneqq}\mleft( s \mright) \) by
\[
  \xi_{\mathfrak{B}}^{:=}\mleft( s \mright)
  \coloneqq s^{m} \mleft( 1 - s \mright)^{m} \xi_{\mathfrak{B}}^{=}\mleft( s \mright),
\]
which does not have the pole of order~\( m \) at~\( s = 0 \),
and still satisfies the functional equation.
By \eqref{xi}~and~\eqref{eqxi}, we have
\[
  \xi_{\mathfrak{B}}^{:=}\mleft( s \mright)
  = \xi_{\mathfrak{B}}^{:}\mleft( s \mright) + \xi_{\mathfrak{B}}^{:}\mleft( 1 - s \mright),
\]
where \( \xi_{\mathfrak{B}}^{:}\mleft( s \mright) \coloneqq
s^{m} \mleft( 1 - s \mright)^{m} \xi_{\mathfrak{B}}\mleft( s \mright) \).
We also define
\( \zeta_{\mathfrak{B}}^{:}\mleft( s \mright) \coloneqq
s^{m} \mleft( 1 - s \mright)^{m} \zeta_{\mathfrak{B}}\mleft( s \mright) \).
In addition to the simple poles given by~\eqref{pPole} when \( n \ne 0 \),
this function has other simple poles from
the \( \Gamma\mleft( s/2 \mright) \) factor in~\eqref{xi}
located at \( -2n \) for~\( n \in \mathbb{N} \),
each with residue~\( -2 \mleft( -1 \mright)^{n} / \Gamma\mleft( n \mright) \).
So the principal part of~\( \xi_{\mathfrak{B}}^{:}\mleft( s \mright) \) is
\begin{equation}
  \label{xiPP}
  \xi_{\mathfrak{B}}^{\colon \text{pp}}\mleft( s \mright)
  = \sum_{j = 1}^{m}\sum_{n = -\infty}^{\infty}
  \frac{\frac{2 n \pi i}{\log p_{j}}
  \mleft( 1 - \frac{2 n \pi i}{\log p_{j}} \mright)
  \xi_{\mathfrak{B} \setminus \mleft\{ p_{j} \mright\}}^{:}
  \mleft( \frac{2 n \pi i}{\log p_{j}} \mright)}
  {\log p_{j} \mleft( s - \frac{2 n \pi i}{\log p_{j}} \mright)}
  + \sum_{n = 1}^{\infty} \mleft( -1 \mright)^{n}
  \frac{2 \pi^{n} \mleft( 2n + 1 \mright) \zeta_{\mathfrak{B}}^{:}\mleft( -2n \mright)}
  {\Gamma\mleft( n \mright) \mleft( s + 2n \mright)}.
\end{equation}
Now we can remove this part to get an entire function
\[
  \xi_{\mathfrak{B}}^{\colon \text{reg}} \coloneqq
  \xi_{\mathfrak{B}}^{:}\mleft( s \mright)
  - \xi_{\mathfrak{B}}^{\colon \text{pp}}\mleft( s \mright).
\]
Finally our approximation to \( \xi\mleft( s \mright) \)~and~\( \zeta\mleft( s \mright) \)
are given by
\begin{equation}
  \label{approx}
  \xi_{\mathfrak{B}}^{\approx}\mleft( s \mright) \coloneqq
  \frac{\xi_{\mathfrak{B}}^{\colon \text{reg}}\mleft( s \mright)
  + \xi_{\mathfrak{B}}^{\colon \text{reg}}\mleft( 1 - s \mright)}
  {s^{m} \mleft( 1 - s \mright)^{m}}
  \text{~and~}
  \zeta_{\mathfrak{B}}^{\approx}\mleft( s \mright) \coloneqq
  \frac{\xi_{\mathfrak{B}}^{\approx}\mleft( s \mright)}{g\mleft( s \mright)},
\end{equation}
where \( g\mleft( s \mright) \) is the factor from~\eqref{xi}.
If \( \mathfrak{B} = \mleft\{ p \colon p \le u \mright\} \) for some \( u \in \mathbb{N} \),
we abbreviate the subscript as \( \xi_{u}^{\approx}\mleft( s \mright), \zeta_{u}^{\approx}\mleft( s \mright) \), etc.

\section{The convergence problem for~
\texorpdfstring{\( \xi_{\mathfrak{B}}^{\colon \text{pp}}\mleft( s \mright) \)}{}}

In the definition of~\( \xi_{\mathfrak{B}}^{\colon \text{pp}}\mleft( s \mright) \) by~\eqref{xiPP},
the second sum over negative even integers is convergent due to the \( \Gamma\mleft( n \mright) \)
in the denominator. For the first sum, we know
\[
  \abs{ \xi_{\mathfrak{B} \setminus \mleft\{ p \mright\}}^{\colon}\mleft( \frac{2 n \pi i}{\log p} \mright) }
  = \frac{1}{2} \abs{ \frac{2n \pi i}{\log p} }^{m}
  \abs{ 1 - \frac{2 n \pi i}{\log p} }^{m}
  \abs{ \Gamma\mleft( \frac{n \pi i}{\log p} \mright) }
  \abs{ \zeta_{\mathfrak{B} \setminus \mleft\{ p \mright\}}\mleft( \frac{2 n \pi i}{\log p} \mright) },
\]
thus as long as the last \( \zeta_{\mathfrak{B} \setminus \mleft\{ p \mright\}} \) factor
is not large, it will also be convergent
due to exponential decay of the \( \Gamma \)-function on the imaginary axis.
However by~\eqref{zeta}, we know
\[
  \abs{ \zeta_{\mathfrak{B} \setminus \mleft\{ p \mright\}}\mleft( \frac{2 n \pi i}{\log p} \mright) }
  = \prod_{\substack{q \in \mathfrak{B} \\ q \ne p}}
  \abs{ 1 - q^{- \frac{2 n \pi i}{\log p}} }^{-1}
  = \prod_{\substack{q \in \mathfrak{B} \\ q \ne p}}
  \abs{ 1 - e^{- 2 \pi i\frac{n \log q}{\log p}} }^{-1},
\]
which could be large when
\begin{equation}
  \label{linLog}
  \frac{n \log q}{\log p} \approx m
  \Leftrightarrow
  m \log p - n \log q = \pm \varepsilon,
\end{equation}
for some small \( \varepsilon > 0 \)~and~\( m, n \in \mathbb{Z} \).
 {Bounds on
linear combinations of logarithms
of natural numbers is a classical problem in Number Theory;
here we use results from~\cite{gouillonExplicitLowerBounds2006}.
In notation from this paper we}
setup \( \alpha_{1} = p, \alpha_{2} = q \)~and~\( b_{1} = m, b_{2} = n \),
whence the degree~\( d \) of~\( \alpha_{1} \) over~\( \mathbb{Q} \) is simply \( 1 \),
and \( a = 1, h\mleft( \alpha_{1} \mright) = \log p, h\mleft( \alpha_{2} \mright) = \log q, D = 1 \)
for similar reasons. For more details during the setup process,
see~\cite[Section~2]{gouillonExplicitLowerBounds2006}.
Next since primes are multiplicatively independent,
we can apply~\cite[Corollary~2.3]{gouillonExplicitLowerBounds2006}
with \( A_{1} = p, A_{2} = q \) and
\begin{align*}
  b
  & = \frac{n}{\log p} + \frac{m}{\log q}
  \approx \frac{2 n}{\log p}, \\
  h
  & = \max \mleft\{ \log b + 3.1, 1000 \mright\},
\end{align*}
which gives the lower bound
\begin{equation}
  \label{epsBound}
  \abs{ \varepsilon }
  > \exp\mleft( -36821 \cdot h \log p \log q \mright).
\end{equation}
Now for this \( q \) in~\eqref{linLog}, its Euler factor is
\begin{align*}
  \abs{ 1 - \exp\mleft( - 2 \pi i \frac{n \log q}{\log p} \mright) }
  & = \abs{ 1 - \exp\mleft( - 2 \pi i \frac{\varepsilon}{\log p} \mright) } \\
  & = \abs{ \frac{2 \pi i \varepsilon}{\log p} - \mleft( \frac{2 \pi i \varepsilon}{\log p} \mright)^{2} / 2! + \ldots } \\
  & > \frac{\pi \varepsilon}{\log p}.
\end{align*}
Then by~\eqref{epsBound} we see
\[
  \abs{ 1 - \exp\mleft( - 2 \pi i \frac{n \log q}{\log p} \mright) }^{-1}
  < \frac{\log p}{\pi} \exp\mleft( 36821 \cdot h \log p \log q \mright).
\]
This shows the convergence of~\eqref{xiPP}, yet the bound is too large for our purpose.

Now we show there are actually cancellations in the first sum of~\eqref{xiPP}.
First rewrite~\eqref{linLog} as
\[
  \frac{n}{\log p} \approx \tau \approx \frac{m}{\log q},
\]
for \( \tau \in \mathbb{R} \), whence we can write
\[
  \frac{n}{\log p} = \tau + \delta_{p} \varepsilon,
  \text{ and }
  \frac{m}{\log q} = \tau + \delta_{q} \varepsilon
\]
for some small \( \varepsilon > 0 \)~and~\( \delta_{p}, \delta_{q} \in \mathbb{R} \).
Now
\begin{align*}
  1 - \exp\mleft( - 2 \pi i \frac{n \log q}{\log p} \mright)
  & = 1 - \exp\mleft( - 2 \pi i \log q \mleft( \tau + \delta_{p} \varepsilon \mright) \mright) \\
  & = 1 - \exp\mleft( - 2 \pi i \log q \mleft( \delta_{p} - \delta_{q} \mright) \varepsilon \mright) \\
  & = 2 \pi i \log q \mleft( \delta_{p} - \delta_{q} \mright) \varepsilon
  + O\mleft( \varepsilon^{2} \mright),
\end{align*}
which along with the \( \log p \) in the denominator of the first sum in~\eqref{xiPP}
gives a term containing
\begin{equation}
  \label{pTerm}
  \frac{1}{2 \pi i \log p \log q} \mleft( \delta_{p} - \delta_{q} \mright)
  \mleft( \frac{1}{\varepsilon} + O\mleft( 1 \mright) \mright).
\end{equation}
Then if we repeat this analysis for~\( q \), the corresponding term contains
\begin{equation}
  \label{qTerm}
  \frac{1}{2 \pi i \log q \log p} \mleft( \delta_{q} - \delta_{p} \mright)
  \mleft( \frac{1}{\varepsilon} + O\mleft( 1 \mright) \mright).
\end{equation}
We see that the main parts (that are due to~\( \varepsilon^{-1} \)) of
\eqref{pTerm}~and~\eqref{qTerm} are of opposite signs,
so taken together,
\eqref{pTerm}~and~\eqref{qTerm} contribute to~\eqref{xiPP}
much smaller amount than  when considered individually.

This argument shows that in~\eqref{xiPP},
the sum defining \( \xi_{\mathfrak{B}}^{\colon \text{pp}}\mleft( s \mright) \)
should not be taken separately but considered as a whole,
since the algebraic structure between the summands is quite intertwined.
This approach is what we will use in the next section.
Such a treatment can also be found in~\cite[Theorem~14.27]{titchmarshTheoryRiemannZetafunction1986},
where the convergence is established along some convenient sequence~\( \mleft( T_{k} \mright) \)
instead of the general~\( T \).

\section{Proof of \autoref{thm:sieveAtOneHalf}~and~\autoref{thm:sieve}}

We need the following two bounds for the \( \Gamma \)-function.

\begin{lemma}
  For \( x > 0 \),
  \[
    \sqrt{2 \pi} x^{x - 1/2} e^{-x}
    \le \Gamma\mleft( x \mright)
    \le \sqrt{2 \pi} x^{x - 1/2} e^{-x} e^{1 / \mleft( 12 x \mright)}.
  \]
\end{lemma}

\begin{proof}
  This is Stirling's series truncated at~\( x^{-1} \).
\end{proof}

\begin{lemma}
  \label{gamma}
  Let \( x, y \in \mathbb{R} \)~and~\( x \ge 1 \), then
  \[
    \abs{ \frac{\Gamma\mleft( x + i y \mright)}{\Gamma\mleft( x \mright)} }
    \le \exp\mleft( -\tfrac{1}{4} \min \mleft\{ \abs{ y }, y^{2}/x \mright\} \mright).
  \]
\end{lemma}

\begin{proof}
  We have the following identity for the \( \Gamma \)-function
  \[
    \abs{ \frac{\Gamma\mleft( x + iy \mright)}{\Gamma\mleft( x \mright)} }
    = \prod_{k = 0}^{\infty} \mleft( 1 + \frac{y^{2}}{\mleft( x + k \mright)^{2}} \mright)^{-1/2}.
  \]
  After taking logarithm we see a lower bound is needed for
  \[
    \frac{1}{2} \sum_{k = 0}^{\infty} \log\mleft( 1 + \frac{y^{2}}{\mleft( x + k \mright)^{2}} \mright).
  \]
  By the elementary inequality~\( \log\mleft( 1 + t \mright) \ge t / \mleft( 1 + t \mright) \) for \( t > -1 \),
  we see
  \[
    \sum_{k = 0}^{\infty} \log\mleft( 1 + \frac{y^{2}}{\mleft( x + k \mright)^{2}} \mright)
    \ge \sum_{k = 0}^{\infty} \frac{y^{2}}{y^{2} + \mleft( x + k \mright)^{2}}.
  \]
  When \( x \ge \abs{ y } \),
  \[
    \sum_{k = 0}^{\infty} \frac{y^{2}}{y^{2} + \mleft( x + k \mright)^{2}}
    \ge \frac{y^{2}}{2} \sum_{k = 0}^{\infty} \frac{1}{\mleft( x + k \mright)^{2}}
    >  \frac{y^{2}}{2} \sum_{k = 0}^{\infty} \frac{1}{\mleft( x + k \mright)\mleft( x + k + 1 \mright)}
    = \frac{y^{2}}{2 x}.
  \]
  When \( x < \abs{ y } \),
  \[
    \sum_{k = 0}^{\infty} \frac{y^{2}}{y^{2} + \mleft( x + k \mright)^{2}}
    \ge \frac{y^{2}}{2} \sum_{k = 0}^{\infty} \frac{1}{\mleft( \abs{ y } + k \mright)^{2}}
    > \frac{\abs{ y }}{2}.
  \]
\end{proof}

We also use some results from~\cite{nastasescuClassApproximationsRiemann2022}.
Let \( \mathcal{C} \) be the rectangular contour given by
\( -\tau \pm iT \)~and~\( 1 + \tau \pm iT \),
which is symmetric with respect to~\( \Re s = 1/2 \).
We write \( \mathcal{C} = \mathcal{L} + \mathcal{R} \),
where the split is on the \( \Re s = 1/2 \) line.
Since \( \xi_{u}^{\colon \text{reg}}\mleft( s \mright) \) is entire we know
\[
  \xi_{u}^{\colon \text{reg}}\mleft( \frac{1}{2} \mright)
  = \frac{1}{2 \pi i} \oint_{\mathcal{C}}
  \frac{\xi_{u}^{\colon \text{reg}}\mleft( s \mright)}{s - 1/2} \df s.
\]
Then by~\( s^{m} \mleft( 1 - s \mright)^{m} \xi_{u}^{\approx}\mleft( s \mright)
= \xi_{u}^{\colon \text{reg}}\mleft( s \mright) + \xi_{u}^{\colon \text{reg}}\mleft( 1 - s \mright) \),
we have
\[
  \xi_{u}^{\approx}\mleft( \frac{1}{2} \mright)
  = 2^{2m + 1} \cdot \frac{1}{2 \pi i} \oint_{\mathcal{C}}
  \frac{\xi_{u}^{\colon \text{reg}}\mleft( s \mright)}{s - 1/2} \df s.
\]
Recall \( \xi_{u}^{\colon \text{reg}}\mleft( s \mright)
= \xi_{u}^{:}\mleft( s \mright) - \xi_{u}^{\colon \text{pp}}\mleft( s \mright)
= s^{m} \mleft( 1 - s \mright)^{m} \xi_{u}\mleft( s \mright) - \xi_{u}^{\colon \text{pp}}\mleft( s \mright)
\), now
\[
  \xi_{u}^{\approx}\mleft( \frac{1}{2} \mright)
  = 2^{2m + 1} \cdot \frac{1}{2 \pi i} \oint_{\mathcal{C}}
  \frac{s^{m} \mleft( 1 - s \mright)^{m} \xi_{u}\mleft( s \mright)}{s - 1/2} \df s
  - 2^{2m + 1} \cdot \frac{1}{2 \pi i} \oint_{\mathcal{C}}
  \frac{\xi_{u}^{\colon \text{pp}}\mleft( s \mright)}{s - 1/2} \df s.
  \numberthis \label{xiAppAtHalf}
\]
We will show the second integral is small.
By the same argument of~\cite[(4.12) \& (4.18)]{nastasescuClassApproximationsRiemann2022},
we know
\begin{equation}
  \label{ppInt}
  \frac{1}{2 \pi i} \oint_{\mathcal{C}}
  \frac{\xi_{u}^{\colon \text{pp}}\mleft( s \mright)}{s - 1/2} \df s
  = \sum_{j = 1}^{m} \sum_{\abs{ n } > T \log p_{j} / \mleft( 2 \pi \mright)}
  \res\mleft( \frac{2 n \pi i}{\log p_{j}} \mright)
  + \sum_{n > \tau/2} \res\mleft( -2n \mright),
\end{equation}
where the poles are those outside~\( \mathcal{C} \),
whose residues are given in~\eqref{xiPP}.
For the \( -2n \) type of poles, their residues are
\[
  \sum_{n > \tau/2} \mleft( -1 \mright)^{n} \frac{2 \pi^{n} \mleft( -2n \mright)^{m} \mleft( 2n + 1 \mright)^{m+1}
  \zeta_{\mathfrak{B}}\mleft( -2n \mright)}{\Gamma\mleft( n \mright)},
\]
which is convergent, hence it tends to~\( 0 \) as we let \( \tau \to \infty \).

Next we deal with the \( 2 n \pi i / \log p \)~type poles.
There exist \( T_{1} \) between \( T+1, T+2 \) and \( \tau > 0 \) such that
\begin{equation}
  \label{poleDist}
  \min_{1 \le j \le m} \min_{n \in \mathbb{Z}}
  \abs{ T_{1} - \frac{2 \pi n}{\log p_{j}} } \ge \frac{2 \pi}{u},
  \text{~and~}
  \min_{n \in \mathbb{Z}} \abs{ \tau + 2n } \ge 1.
\end{equation}
This is so because in an
fixed interval of length~\( 1 \),
the number of poles is
\[
  \sum_{j = 1}^{m} \frac{\log p_{j}}{2 \pi} \le \frac{u}{2 \pi}.
\]
Let \( \mathcal{C}_{1} \) be the rectangular contour given by
\( - \tau \pm i T_{1} \) and \( 1 + \tau \pm i T_{1} \).
By the residue theorem we see for those poles between \( \mathcal{C} \)~and~\( \mathcal{C}_{1} \),
\[
  \sum_{j = 1}^{m} \sum_{T < \frac{2 \pi n}{\log p_{j}} < T_{1}}
  \res\mleft( \frac{2 n \pi i}{\log p_{j}} \mright)
  = \frac{1}{2 \pi i} \oint_{\mathcal{C}_{1}}
  \frac{s^{m} \mleft( 1 - s \mright)^{m} \xi_{u}\mleft( s \mright)}{s - 1/2} \df s
  - \frac{1}{2 \pi i} \oint_{\mathcal{C}}
  \frac{s^{m} \mleft( 1 - s \mright)^{m} \xi_{u}\mleft( s \mright)}{s - 1/2} \df s.
\]
The integral over~\( \mathcal{C}_{1} \) is
\begin{equation}
  \label{iOne}
  I_{1} \le \frac{\abs{ \mathcal{C}_{1} }}{2 \pi}
  \frac{\mleft( 1 + \tau^{2} + T_{1}^{2} \mright)^{m + 1}}{T_{1}}
  \pi^{\frac{1 + \tau}{2}} \Gamma\mleft( \tfrac{1 + \tau}{2} \mright)
  \exp\mleft( -\frac{1}{4} \min \mleft\{ \abs{ T_{1} }, \frac{T_{1}^{2}}{1 + \tau} \mright\} \mright)
  \abs{ \zeta_{u}\mleft( i T_{1} \mright) },
\end{equation}
where
\[
  \abs{ \zeta_{u}\mleft( i T_{1} \mright) }
  = \prod_{j = 1}^{m} \frac{1}{\abs{ 1 - e^{- iT_{1} \log p_{j}} }}
  = \prod_{j = 1}^{m} \frac{1}{\abs{ 2 \sin\mleft( \frac{T_{1} \log p_{j}}{2} \mright) }}
  \le \mleft( \frac{\pi u}{4} \mright)^{m} \prod_{j = 1}^{m} \frac{1}{\log p_{j}}
\]
by~\eqref{poleDist} and the following inequality
\[
  \frac{1}{\abs{ \sin\mleft( \pi x \mright) }}
  \le \frac{\pi}{2 \norm{x}},
\]
where \( \norm{ x } \) is the distance to the nearest integer.
Then this argument can be repeated to find \( T_{2}, T_{3}, \ldots \)
such that \( T_{k} \to \infty \) as~\( k \to \infty \) and
\[
  \sum_{j = 1}^{m} \sum_{T < \frac{2 \pi n}{\log p_{j}} < T_{k}}
  \res\mleft( \frac{2 n \pi i}{\log p_{j}} \mright)
  = \sum_{j = 1}^{k} I_{j},
\]
where \( I_{j} \) is the corresponding bound~\eqref{iOne} for each integral over~\( \mathcal{C}_{j} \).
By the exponential decay in the height~\( T_{j} \) in~\eqref{iOne}, the sum~\( \sum_{j} I_{j} \)
is convergent. Moreover, if we choose the height~\( T \) of
the initial contour~\( \mathcal{C} \) to satisfy~\eqref{poleDist} as well,
a similar bound as~\eqref{iOne} will apply to the integral over~\( \mathcal{C} \) too.
Thus we conclude the principal part integral~\eqref{ppInt} tends to \( 0 \) as we let \( T \to \infty \).

Since \( s^{m} \mleft( 1 - s \mright)^{m} \xi\mleft( s \mright) \) is also entire
and symmetric to~\( \Re s = 1/2 \), we know
\begin{align*}
  \xi\mleft( \frac{1}{2} \mright)
  & = 2^{2m} \cdot \frac{1}{2 \pi i} \oint_{\mathcal{C}}
  \frac{s^{m} \mleft( 1 - s \mright)^{m} \xi\mleft( s \mright)}{s - 1/2} \df s \\
  & = 2^{2m + 1} \cdot \frac{1}{2 \pi i} \int_{\mathcal{R}}
  \frac{s^{m} \mleft( 1 - s \mright)^{m} \xi\mleft( s \mright)}{s - 1/2} \df s.
  \numberthis \label{xiAtHalf}
\end{align*}
Taking the difference between the first integral in~\eqref{xiAppAtHalf}~
and~\eqref{xiAtHalf} we obtain
\begin{multline}
  \label{xidfLR}
  \xi\mleft( \frac{1}{2} \mright)
  - \xi_{u}^{\approx}\mleft( \frac{1}{2} \mright)
  = - 2^{2m + 1} \cdot \frac{1}{2 \pi i} \int_{\mathcal{L}}
  s^{m} \mleft( 1 - s \mright)^{m}
  \frac{\xi_{u}\mleft( s \mright)}{s - 1/2} \df s \\
  + 2^{2m + 1} \cdot \frac{1}{2 \pi i} \int_{\mathcal{R}}
  s^{m} \mleft( 1 - s \mright)^{m}
  \frac{\xi\mleft( s \mright) - \xi_{u}\mleft( s \mright)}{s - 1/2} \df s.
\end{multline}

The contribution on the horizontal part of~\( \mathcal{L} \) is
small due to the \( \Gamma \)-function,
so the integral over~\( \mathcal{L} \) is actually
\begin{multline}
  \label{leftInt}
  2^{2m + 1} \cdot \frac{1}{2 \pi i} \int_{\mleft( - \tau \mright)}
  s^{m} \mleft( 1 - s \mright)^{m}
  \frac{\xi_{u}\mleft( s \mright)}{s - 1/2} \df s \le 2^{2m} \frac{1}{2 \pi i}
  \int_{-\infty}^{\infty} \abs{ \mleft( - \tau + it \mright) \mleft( 1 + \tau - it \mright) }^{m + 1} \\
  \times \pi^{- \tau/2} \abs{ \Gamma\mleft( \frac{- \tau + it}{2} \mright) }
  \abs{ \frac{\zeta_{u}\mleft( - \tau + it \mright)}
  {\tau + 1/2 - it} } \df t.
\end{multline}
Trivially bounding \( \zeta_{u}\mleft( - \tau + it \mright) \) by
\[
  \abs{ \zeta_{u}\mleft( - \tau  + it \mright) }
  = \prod_{p \le u} \abs{ \frac{1}{1 - p^{\tau - it}} }
  \le \prod_{p \le u} \abs{ \frac{1}{p^{\tau} - 1} },
\]
using Stirling's approximation to the \( \Gamma \)-function
\[
  \abs{ \Gamma\mleft( - \tau + it \mright) }
  = \sqrt{2 \pi} \abs{ t }^{- \tau - 1/2} e^{- \pi \abs{ t }/2}
  \mleft( 1 + O\mleft( 1/ \abs{ t } \mright) \mright),
\]
we see the function in~\eqref{leftInt} is integrable in~\( t \)
and the integral goes to~\( 0 \) as~\( - \tau \to - \infty \) uniformly.
Now we are left with the integral over~\( \mathcal{R} \) in~\eqref{xidfLR},
and similarly the honrizontal contribution can be ignored,
\[
  2^{2m + 1} \cdot \frac{1}{2 \pi i} \int_{\mleft( \sigma \mright)}
  s^{m} \mleft( 1 - s \mright)^{m} g\mleft( s \mright)
  \frac{\zeta\mleft( s \mright) - \zeta_{u}\mleft( s \mright)}{s - 1/2} \df s.
\]
Let \( P^{-}\mleft(n\mright) \) denote the largest prime factor of~\( n \)
(with \( P^{-}\mleft(1\mright) = 1 \)). Since
\begin{equation}
  \label{seriesdf}
  \zeta\mleft( s \mright) - \zeta_{u}\mleft( s \mright)
  = \sum_{P^{-}\mleft( n \mright) > u} n^{-s},
\end{equation}
when \( \sigma > 1 \) we can exchange integration with summation to get
\begin{equation}
  \label{sumJn}
  \mleft( -1 \mright)^{m} 2^{2m + 1} \cdot \frac{1}{2 \pi i} \int_{\mleft( \sigma \mright)}
  s^{m + 1} \mleft( s - 1 \mright)^{m + 1} g\mleft( s \mright)
  \frac{\zeta\mleft( s \mright) - \zeta_{u}\mleft( s \mright)}{s - 1/2} \df s
  = \mleft( -1 \mright)^{m} 2^{2m + 1} \sum_{P^{-}\mleft( n \mright) > u} J\mleft( n \mright)
\end{equation}
where
\begin{equation}
  \label{jInt}
  J\mleft( n \mright) =
  \frac{1}{2 \pi i} \int_{\mleft( \sigma \mright)}
  \frac{s^{m + 1} \mleft( s - 1 \mright)^{m + 1}}{s - 1/2}
  \Gamma\mleft( s/2 \mright) \mleft( \pi n^{2} \mright)^{-s/2} \df s.
\end{equation}

Next we split the line integral in~\eqref{jInt} into two parts \( J_{1} \)~and~\( J_{2} \),
where \( J_{1} \) is taken over~\( \abs{ \Im s } \le t_{0} \) for some~\( t_{0} > 0 \)
to be chosen later, and \( J_{2} \) the rest.

We choose \( \sigma = \sigma_{n} = 2 \pi n^{2} \)~and~\( t_{0} = t_{n} = n \log^{3/4} n \)
on~\( J_{1} \), whence
\[
  \abs{ s - 1 }^{2m + 1} < \abs{ \frac{s^{m+1} \mleft( s - 1 \mright)^{m+1}}{s - 1/2} }
  < \abs{ s }^{2m + 1}.
\]
By the elementary inequality~\( e^{x} \ge \mleft( 1 + x/n \mright)^{n} \ge 1 + x \),
for upper bound we have
\begin{equation}
  \label{sUpper}
  \abs{ s }^{2m + 1}
  < \sigma^{2m + 1} \mleft( 1 + \frac{t_{0}}{\sigma} \mright)^{2 m + 1}
  \le \sigma^{2m + 1} \mleft( 1 + \frac{\log^{3/4} n}{2 \pi n} \mright)
  \exp\mleft( \frac{m \log^{3/4} n}{\pi n} \mright),
\end{equation}
and for lower bound
\begin{equation}
  \label{sLower}
  \abs{ s - 1 }^{2m + 1}
  > \sigma^{2m + 1} \mleft( 1 - \frac{1}{\sigma} \mright)^{2m + 1}
  \ge \sigma^{2m + 1} \mleft( 1 - \frac{1}{2 \pi n^{2}} \mright)
  \mleft( 1 - \frac{m}{ \pi n^{2}} \mright).
\end{equation}

Now we can replace the rational part of~\( s \) in~\eqref{jInt}
by~\( \sigma^{2m + 1} \) with the two bounds given by \eqref{sUpper}~and~\eqref{sLower}.
We first integrate along the whole \( \Re s = \sigma \)~line to get
\begin{equation}
  \label{oneMain}
  \frac{1}{2 \pi i} \int_{\mleft( \sigma \mright)} \Gamma\mleft( s/2 \mright)
  \mleft( \pi n^{2} \mright)^{-s/2} \df s
  = 2 e^{- \pi n^{2}},
\end{equation}
then subtract the \( \abs{ t } > t_{0} \)~part, which is
\begin{multline}
  \label{tError}
  \frac{1}{2 \pi} \int_{\abs{ t } > t_{0}}
  \Gamma\mleft( \sigma/2 + it/2 \mright) \mleft( \pi n^{2} \mright)^{-\sigma/2 - it/2} \df t \\
  \le \frac{1}{\pi} \Gamma\mleft( \sigma/2 \mright) \mleft( \pi n^{2} \mright)^{- \sigma/2}
  \mleft( \int_{t_{0}}^{\sigma} e^{-t^{2}/\mleft( 4 \sigma \mright)} \df t
  + \int_{\sigma}^{\infty} e^{-t/4} \df t \mright) \\
  \le \frac{1}{\pi} \Gamma\mleft( \sigma/2 \mright) \mleft( \pi n^{2} \mright)^{- \sigma/2}
  \mleft( \frac{4 \pi n}{\log^{3/4} n} e^{- \frac{\log^{3/2} n}{8 \pi}}
  + 4 e^{- \pi n^{2} / 2} \mright),
\end{multline}
where the bound for the integral from \( t_{0} \) to~\( \sigma \) comes from
\[
  \int_{t_{0}}^{\sigma} e^{-t^{2}/\mleft( 4 \sigma \mright)} \df t
  < 2 \sqrt{\sigma} \int_{\frac{t_{0}}{2 \sqrt{\sigma}}}^{\infty}
  e^{-t^{2}} \df t
  < 2 \sqrt{\sigma} \int_{\frac{t_{0}}{2 \sqrt{\sigma}}}^{\infty}
  \frac{t}{t_{0}/\mleft( 2 \sqrt{\sigma} \mright)} e^{-t^{2}} \df t
  = \frac{2 \sigma}{t_{0}} e^{-t_{0}^{2} / \mleft( 4 \sigma \mright)}.
\]
By~\eqref{sUpper}, \eqref{oneMain}~and~\eqref{tError},
the upper bound for~\( J_{1} \) is
\begin{multline}
  \label{oneUpperBound}
  \frac{\abs{ J_{1} }}{2 e^{- \pi n^{2}} \sigma^{2m + 1}}
  \le \mleft( 1 + 2 \sqrt{2} e^{\frac{1}{12 \pi n^{2}}}
  \mleft( \frac{1}{\log^{3/4} n} e^{- \frac{\log^{3/2} n}{8 \pi}}
  + \frac{1}{\pi n} e^{- \pi n^{2} / 2} \mright) \mright) \\
  \times \mleft( 1 + \frac{\log^{3/4} n}{2 \pi n} \mright)
  \exp\mleft( \frac{m \log^{3/4} n}{\pi n} \mright),
\end{multline}
and the lower is
\begin{multline}
  \label{oneLowerBound}
  \frac{\abs{ J_{1} }}{2 e^{- \pi n^{2}} \sigma^{2m + 1}}
  \ge \mleft( 1 - 2 \sqrt{2} e^{\frac{1}{12 \pi n^{2}}}
  \mleft( \frac{1}{\log^{3/4} n} e^{- \frac{\log^{3/2} n}{8 \pi}}
  + \frac{1}{\pi n} e^{- \pi n^{2} / 2} \mright) \mright) \\
  \times \bigg( 1 - \frac{1}{2 \pi n^{2}} \bigg)
  \bigg( 1 - \frac{m}{ \pi n^{2}} \bigg),
\end{multline}
by~\eqref{sLower}.

On~\( J_{2} \) we need the following bound
\[
  \abs{ \frac{s^{m + 1} \mleft( s - 1 \mright)^{m+1}}{s - 1/2} }
  < 2 \abs{ s }^{2m + 1}
  < 2 \sigma^{2m + 1}
  \exp\mleft( \mleft( 2m + 1 \mright) \log \mleft( 1 + \frac{t}{\sigma} \mright) \mright)
  < 2 \sigma^{2m + 1} e^{\mleft( 2m + 1 \mright) t / \sigma},
\]
then we can bound \( J_{2} \) by
\begin{align*}
  \abs{ J_{2} }
  & \le \frac{2}{\pi} \mleft( \pi n^{2} \mright)^{-\sigma/2}
  \Gamma\mleft( \sigma/2 \mright) \sigma^{2m + 1}
  \int_{t_{0}}^{\infty} e^{\mleft( 2m + 1 \mright)t/\sigma
  - \frac{1}{4} \min\mleft\{ t, t^{2}/\sigma \mright\}} \df t \\
  & \le \frac{2 \sqrt{2} \sigma^{2m + 1}}{\pi n} e^{- \pi n^{2} + 1/\mleft( 12 \pi n^{2} \mright)}
  \mleft( \int_{t_{0}}^{\sigma} e^{\mleft( 2m + 1 \mright)t / \sigma - t^{2}/\mleft( 4 \sigma \mright)}
  \df t + \int_{\sigma}^{\infty} e^{\mleft( 2m + 1 \mright)t/\sigma - t/4} \df t \mright),
\end{align*}
where the first integral is
\begin{multline*}
  e^{\frac{\mleft( 2m + 1 \mright)^{2}}{\sigma}}
  \int_{t_{0} + 2 \mleft( 2m + 1 \mright)}^{\infty} e^{-t^{2} / 4 \sigma} \df t
  = 2 \sqrt{\sigma} e^{\frac{\mleft( 2m + 1 \mright)^{2}}{\sigma}}
  \int_{\frac{t_{0} + 2 \mleft( 2m + 1 \mright)}{2 \sqrt{\sigma}}}^{\infty} e^{-t^{2}} \df t \\
  < \frac{2 \sigma}{t_{0} + 2 \mleft( 2m + 1 \mright)}
  \exp\mleft(- \frac{t_{0}^{2} + 2\mleft( 2m + 1 \mright)t_{0}}{4 \sigma} \mright) \\
  = \frac{4 \pi n^{2}}{n \log^{3/4} n + 2\mleft( 2m + 1 \mright)}
  \exp\mleft( -\frac{n \log^{3/2}n + 2\mleft( 2m + 1 \mright)\log n}{8 \pi n} \mright),
\end{multline*}
and the second is
\[
  \frac{4 \sigma}{\sigma - 4 \mleft( 2m + 1 \mright)} e^{2m + 1 - \sigma/4}
  = \frac{8 \pi n^{2}}{2 \pi n^{2} - 4 \mleft( 2m + 1 \mright)}
  \exp\mleft( 2m + 1 - \frac{\pi n^{2}}{2} \mright).
\]
So together
\begin{multline}
  \label{twoBound}
  \frac{\abs{ J_{2} }}{2 e^{-\pi n^{2}} \sigma^{2m + 1}}
  \le \frac{\sqrt{2}}{\pi n} e^{\frac{1}{12 \pi n^{2}}}
  \mleft( \frac{4 \pi n^{2}}{n \log^{3/4} n + 2\mleft( 2m + 1 \mright)} \mright. \\
  \times \mleft. \exp\mleft( -\frac{n \log^{3/2}n + 2\mleft( 2m + 1 \mright)\log n}{8 \pi n} \mright)
  + \frac{8 \pi n^{2}}{2 \pi n^{2} - 4 \mleft( 2m + 1 \mright)}
  \exp\mleft( 2m + 1 - \frac{\pi n^{2}}{2} \mright) \mright).
\end{multline}

Finally by~\( J\mleft( n \mright) = J_{1} + J_{2} \), we obtain
\begin{equation}
  \label{jBound}
  1 + L_{n} \le \frac{\abs{ J\mleft( n \mright) }}{2e^{-\pi n^{2}}
  \mleft( 2 \pi n^{2} \mright)^{2m + 1}} \le 1 + U_{n},
\end{equation}
where \( 1 + U_{n} \) is given by adding the right side of \eqref{oneUpperBound}~and~\eqref{twoBound},
while \( 1 + L_{n} \) by adding the right side of~\eqref{oneLowerBound} then subtracting that of~\eqref{twoBound},
and both~\( U_{n}, L_{n} \to 0 \) monotonically as~\( n \to \infty \).
Recall \( m \) is the number of primes used in the approximation.
When the series differs from the \( n \)-th~term we know \( m = \pi\mleft( n - 1 \mright) \),
since the primes we used must be smaller than~\( n \) by~\eqref{seriesdf}.

Then by~\eqref{sumJn}
\begin{equation}
  \label{pQuality}
  \abs{ \xi\mleft( 1/2 \mright)
  - \xi_{\mathfrak{B}}^{\approx}\mleft( 1/2 \mright) }
  - 2^{2m + 1} \abs{ J\mleft( p_{m + 1} \mright) }
  \le 2^{2m + 1} \sum_{n > p_{m+1}} \abs{ J\mleft( n \mright) },
\end{equation}
which will be exponentially smaller than~\( 2^{2m+1} \abs{ J\mleft( p_{m+1} \mright) } \)
as shown in~\cite[11]{nastasescuClassApproximationsRiemann2022}.
Then let \( \mathfrak{B} = \mleft\{ 2, \ldots, p_{m} \mright\} \),
and let \( \ell \) be any composite number such that \( p_{m} < \ell < p_{m+1} \).
Since \( P^{+}\mleft( \ell \mright) < u \), we know
\[
  \zeta\mleft( s \mright) - \frac{1}{1 - \ell^{-s}} \zeta_{u}\mleft( s \mright)
  = \sum_{P\mleft( n \mright) > u} n^{-s}
  - \sum_{\substack{\ell \mid m \\ P\mleft( m/\ell \mright) < u}} m^{-s}
  = -\frac{1}{\ell^{s}} + \ldots.
\]
The same argument now gives
\begin{equation}
  \label{kQuality}
  \abs{ \xi\mleft( 1/2 \mright)
  - \xi_{\mathfrak{B} \cup \mleft\{ \ell \mright\}}^{\approx}\mleft( 1/2 \mright) }
  - 2^{2m+1} \abs{ J\mleft( \ell \mright) }
  \le 2^{2m + 1} \sum_{n > \ell} \abs{ J\mleft( \ell \mright) }.
\end{equation}
Comparing \eqref{pQuality} to~\eqref{kQuality} with our choice~\( p_{m} < \ell < p_{m+1} \),
we see \eqref{pQuality} is always smaller.

Next we use the bounds~\eqref{jBound} to prove~\autoref{thm:sieveAtOneHalf}.
Precisely, this requires the upper bounds at primes being smaller than
the lower bounds at composite numbers. However, the lower correction factor~\( 1 + L_{n} \)
is negative up to~\( n = 367 \), which becomes vacuous since absolute values are always nonnegative.
Thus the primes strictly smaller than~\( 367 \) are checked manually.
Relevant bounds are listed in~\autoref{tbl:approx},
and the actual values in~\cite[Table~2]{matiyasevichFewFactorsEuler2017}
lie in between.

\begin{table}[htb]
  \centering
  \caption{Bound for~\( { \xi_{\mathfrak{B}}^{\approx}\mleft( 1/2 \mright) - \xi\mleft( 1/2 \mright) } \)}
  \begin{tabular}{lll}
    \toprule
    \( \mathfrak{B} \) & Lower Bound & Upper Bound \\
    \midrule
    \( \mathfrak{B}_{\le 2} \) & \( -4.65 \ldots \times 10^{-6} \) & \( 8.81 \ldots \times 10^{-6} \) \\
    \( \mathfrak{B}_{\le 2} \cup \mleft\{ 3 \mright\} \) & \( -1.02 \ldots \times 10^{-21} \) & \( 2.32 \ldots \times 10^{-21} \) \\
    \( \mathfrak{B}_{\le 3} \cup \mleft\{ 5 \mright\} \) & \( -1.62 \ldots \times 10^{-47} \) & \( 4.22 \ldots \times 10^{-47} \) \\
    \( \mathfrak{B}_{\le 5} \cup \mleft\{ 6 \mright\} \) & \( -1.93 \ldots \times 10^{-25} \) & \( 5.82 \ldots \times 10^{-25} \) \\
    \( \mathfrak{B}_{\le 5} \cup \mleft\{ 7 \mright\} \) & \( -1.01 \ldots \times 10^{-136} \) & \( 2.89 \ldots \times 10^{-136} \) \\
    \( \mathfrak{B}_{\le 7} \cup \mleft\{ 10 \mright\} \) & \( -1.21 \ldots \times 10^{-102} \) & \( 3.82 \ldots \times 10^{-102} \) \\
    \( \mathfrak{B}_{\le 7} \cup \mleft\{ 11 \mright\} \) & \( -2.65 \ldots \times 10^{-194} \) & \( 8.26 \ldots \times 10^{-194} \) \\
    \( \mathfrak{B}_{\le 11} \cup \mleft\{ 12 \mright\} \) & \( -1.80 \ldots \times 10^{-154} \) & \( 6.20 \ldots \times 10^{-154} \) \\
    \( \mathfrak{B}_{\le 11} \cup \mleft\{ 13 \mright\} \) & \( -2.14 \ldots \times 10^{-348} \) & \( 7.13 \ldots \times 10^{-348} \) \\
    \( \mathfrak{B}_{\le 13} \cup \mleft\{ 14 \mright\} \) & \( -6.16 \ldots \times 10^{-217} \) & \( 2.27 \ldots \times 10^{-216} \) \\
    \( \mathfrak{B}_{\le 13} \cup \mleft\{ 15 \mright\} \) & \( -1.32 \ldots \times 10^{-255} \) & \( 4.80 \ldots \times 10^{-255} \) \\
    \( \mathfrak{B}_{\le 13} \cup \mleft\{ 17 \mright\} \) & \( -4.23 \ldots \times 10^{-438} \) & \( 1.51 \ldots \times 10^{-437} \) \\
    \( \mathfrak{B}_{\le 17} \cup \mleft\{ 18 \mright\} \) & \( -4.00 \ldots \times 10^{-381} \) & \( 1.54 \ldots \times 10^{-380} \) \\
    \( \mathfrak{B}_{\le 17} \cup \mleft\{ 19 \mright\} \) & \( -3.22 \ldots \times 10^{-657} \) & \( 1.21 \ldots \times 10^{-656} \) \\
    \( \mathfrak{B}_{\le 19} \cup \mleft\{ 20 \mright\} \) & \( -6.82 \ldots \times 10^{-476} \) & \( 2.80 \ldots \times 10^{-475} \) \\
    \( \mathfrak{B}_{\le 19} \cup \mleft\{ 21 \mright\} \) & \( -5.00 \ldots \times 10^{-531} \) & \( 2.03 \ldots \times 10^{-530} \) \\
    \( \mathfrak{B}_{\le 19} \cup \mleft\{ 22 \mright\} \) & \( -6.27 \ldots \times 10^{-589} \) & \( 2.53 \ldots \times 10^{-588} \) \\
    \( \mathfrak{B}_{\le 19} \cup \mleft\{ 23 \mright\} \) & \( -1.80 \ldots \times 10^{-1071} \) & \( 7.10 \ldots \times 10^{-1071} \) \\
    \( \mathfrak{B}_{\le 23} \cup \mleft\{ 24 \mright\} \) & \( -2.50 \ldots \times 10^{-705} \) & \( 1.07 \ldots \times 10^{-704} \) \\
    \( \mathfrak{B}_{\le 23} \cup \mleft\{ 26 \mright\} \) & \( -2.63 \ldots \times 10^{-840} \) & \( 1.11 \ldots \times 10^{-839} \) \\
    \( \mathfrak{B}_{\le 23} \cup \mleft\{ 28 \mright\} \) & \( -2.61 \ldots \times 10^{-986} \) & \( 1.09 \ldots \times 10^{-985} \) \\
    \( \mathfrak{B}_{\le 23} \cup \mleft\{ 29 \mright\} \) & \( -5.83 \ldots \times 10^{-1226} \) & \( 2.44 \ldots \times 10^{-1225} \) \\
    \( \mathfrak{B}_{\le 29} \cup \mleft\{ 30 \mright\} \) & \( -3.06 \ldots \times 10^{-1135} \) & \( 1.35 \ldots \times 10^{-1134} \) \\
    \( \mathfrak{B}_{\le 29} \cup \mleft\{ 31 \mright\} \) & \( -5.84 \ldots \times 10^{-1771} \) & \( 2.55 \ldots \times 10^{-1770} \) \\
    \( \mathfrak{B}_{\le 31} \cup \mleft\{ 33 \mright\} \) & \( -5.86 \ldots \times 10^{-1383} \) & \( 2.71 \ldots \times 10^{-1382} \) \\
    \( \mathfrak{B}_{\le 31} \cup \mleft\{ 34 \mright\} \) & \( -1.00 \ldots \times 10^{-1473} \) & \( 4.63 \ldots \times 10^{-1473} \) \\
    \( \mathfrak{B}_{\le 31} \cup \mleft\{ 35 \mright\} \) & \( -3.08 \ldots \times 10^{-1567} \) & \( 1.41 \ldots \times 10^{-1566} \) \\
    \( \mathfrak{B}_{\le 31} \cup \mleft\{ 36 \mright\} \) & \( -1.69 \ldots \times 10^{-1663} \) & \( 7.73 \ldots \times 10^{-1663} \) \\
    \bottomrule
  \end{tabular}
  \label{tbl:approx}
\end{table}

Now \( 367 \)~and~\( 373 \) are the next two consecutive primes,
and the numerical bounds in~\autoref{tbl:367}
demonstrates that the sieve works between them.
\begin{table}[htp]
  \centering
  \caption{First few nontrivial bounds for~\( \abs{ \xi_{\mathfrak{B}}^{\approx}\mleft( 1/2 \mright) - \xi\mleft( 1/2 \mright) } \)}
  \begin{tabular}{lll}
    \toprule
    \( \mathfrak{B} \) & Lower Bound & Upper Bound \\
    \midrule
    \( \mathfrak{B}_{\le 367} \)
    & \( 1.02 \ldots \times 10^{-188909} \) & \( 1.06 \ldots \times 10^{-188906} \) \\
    \( \mathfrak{B}_{\le 367} \cup \mleft\{ 368 \mright\} \)
    & \( 6.88 \ldots \times 10^{-183844} \) & \( 6.04 \ldots \times 10^{-183841} \) \\
    \( \mathfrak{B}_{\le 367} \cup \mleft\{ 369 \mright\}\)
    & \( 4.69 \ldots \times 10^{-184849} \) & \( 3.86 \ldots \times 10^{-184846} \) \\
    \( \mathfrak{B}_{\le 367} \cup \mleft\{ 370 \mright\}\)
    & \( 5.95 \ldots \times 10^{-185857} \) & \( 4.60 \ldots \times 10^{-185854} \) \\
    \( \mathfrak{B}_{\le 367} \cup \mleft\{ 371 \mright\}\)
    & \( 1.40 \ldots \times 10^{-186867} \) & \( 1.02 \ldots \times 10^{-186864} \) \\
    \( \mathfrak{B}_{\le 367} \cup \mleft\{ 372 \mright\}\)
    & \( 6.12 \ldots \times 10^{-187881} \) & \( 4.22 \ldots \times 10^{-187878} \) \\
    \( \mathfrak{B}_{\le 373} \)
    & \( 8.44 \ldots \times 10^{-195051} \) & \( 3.23 \ldots \times 10^{-195048} \) \\
    \bottomrule
  \end{tabular}
  \label{tbl:367}
\end{table}

To show this holds true for larger primes, we observe that the bound~\eqref{jBound}
consists of two parts, an exponentially decaying main term,
multiplied by upper and lower correction factors which monotonically tend to~\( 1 \).
Now we list the value of these two correction factors separately in~\autoref{tbl:ulBound}.
\begin{table}[htp]
  \centering
  \caption{Value of and ratio between \( 1 + L_{n} \)~and~\( 1 + U_{n} \)}
  \begin{tabular}{llll}
    \toprule
    \( n \) & \( 1 + L_{n} \) & \( 1 + U_{n} \)
    & \( \frac{1 + U_{n}}{1 + L_{n}} \) \\
    \midrule
    \( 368 \)
    & \( 1.359 \ldots \times 10^{-3} \) & \( 2.3858 \ldots \) & \( 1755.51 \ldots \) \\
    \( 369 \)
    & \( 1.5455 \ldots \times 10^{-3} \) & \( 2.3845 \ldots \) & \( 1542.81 \ldots \) \\
    \( 370 \)
    & \( 1.7325 \ldots \times 10^{-3} \) & \( 2.3832 \ldots \) & \( 1375.55 \ldots \) \\
    \( 371 \)
    & \( 1.9200 \ldots \times 10^{-3} \) & \( 2.3819 \ldots \) & \( 1240.56 \ldots \) \\
    \( 372 \)
    & \( 2.1080 \ldots \times 10^{-3} \) & \( 2.3806 \ldots \) & \( 1129.35 \ldots \) \\
    \( 373 \)
    & \( 2.2964 \ldots \times 10^{-3} \) & \( 2.3793 \ldots \) & \( 1036.13 \ldots \) \\
    \( 379 \)
    & \( 6.2015 \ldots \times 10^{-3} \) & \( 2.3748 \ldots \) & \( 382.94 \ldots \) \\
    \bottomrule
  \end{tabular}
  \label{tbl:ulBound}
\end{table}
Note by our construction, \( n \) is the last composite number if \( \mathfrak{B} \)
contains one, while \( n \) is the next prime when \( \mathfrak{B} \) does not.
In~\autoref{tbl:ulBound}, we see the ratio between correction factors decreases below~\( 10^{3} \)
after~\( 373 \). Meanwhile in~\autoref{tbl:367} the ratio between
consecutive exponential main terms increases above~\( 10^{1000} \).
The prime counting function~\( \pi\mleft( x \mright) \) in~\( m = \pi\mleft( n - 1 \mright) \)
can be bounded by explicit estimates in \citeauthor{dusartExplicitEstimatesFunctions2018}
~\cite[Corollary~5.2]{dusartExplicitEstimatesFunctions2018}.
By the monotonicity of the bounds~\eqref{jBound}, the sieve will always
pick primes after~\( 373 \). This together with our manual check of
primes smaller than~\( 373 \) establishes~\autoref{thm:sieveAtOneHalf}.

The proof of~\autoref{thm:sieve} is similar, with a general~\( s_{0} \in \mathbb{C} \)
in place of~\( 1/2 \). The estimates above may not work from the beginning,
thus the sieve at this~\( s_{0} \) may need the first \( k_{0} \)~primes to be announced,
for some~\( k_{0} \in \mathbb{N} \) depending on~\( s_{0} \).

\section{Proof of \autoref{thm:limit}~and~\autoref{thm:sign}}

We define
\[
  \xi_{\mleft[ m \mright]}\mleft( s \mright)
  \coloneqq \lim_{\min \mleft\{ p_{1}, \ldots, p_{m} \mright\} \to \infty}
  \xi_{\mleft\{ p_{1}, \ldots, p_{m} \mright\}}^{\approx}\mleft( s \mright),
\]
then by~\eqref{approx} it suffices to show that
\begin{equation}
  \label{xi_m}
  \xi_{\mleft[ m \mright]}\mleft( s_{0} \mright)
  = \frac{1}{s_{0}^{m} \mleft( 1 - s_{0} \mright)^{m}} \cdot
  \frac{1}{2 \pi i} \int_{\mleft( \abs{ \Re s_{0} } + 2 \mright)}
  s^{m} \mleft( 1 - s \mright)^{m}
  g\mleft( s \mright) a\mleft( s, s_{0} \mright) \df s.
\end{equation}
We begin with the following formula~\cite[(4.27)]{nastasescuClassApproximationsRiemann2022},
\begin{equation}
  \label{dfXi}
  \xi\mleft( s_{0} \mright) - \xi_{\mathfrak{B}}^{\approx}\mleft( s_{0} \mright)
  = \frac{1}{s_{0}^{m} \mleft( 1 - s_{0} \mright)^{m}} \cdot
  \frac{1}{2 \pi i} \int_{\mleft( \sigma \mright)}
  s^{m} \mleft( 1 - s \mright)^{m} g\mleft( s \mright)
  \mleft( \zeta\mleft( s \mright) - \zeta_{\mathfrak{B}}\mleft( s \mright) \mright)
  a\mleft( s, s_{0} \mright) \df s,
\end{equation}
where \( \sigma > \max \mleft\{ 1, \abs{ \Re s_{0} }, \abs{ 1 - \Re s_{0} } \mright\} \).
Now \( \mathfrak{B} = \mleft\{ p_{1}, \ldots, p_{m} \mright\} \),
let \( p_{0} \) be the previous prime to~\( p_{1} \).
Then the condition~\( \min \mleft\{ p_{1}, \ldots, p_{m} \mright\} \to \infty \)
is equivalent to~\( p_{0} \to \infty \),
whence the difference in the series on the right hand side of~\eqref{dfXi}
contains
\[
  \lim_{p_{0} \to \infty} \sum_{P^{+}\mleft( n \mright) \le p_{0}} \frac{1}{n^{s}}
  = \lim_{p_{0} \to \infty}
  \frac{1}{2^{s}} + \frac{1}{3^{s}} + \ldots + \frac{1}{p_{0}^{s}} + \ldots
  = \zeta\mleft( s \mright) - 1.
\]
This establishes~\eqref{xi_m} and proves~\autoref{thm:limit}.

In~\cite[Theorem~2]{nastasescuClassApproximationsRiemann2022},
sign alternation at prime numbers is proved, thus it suffices
to show that composite numbers give the opposite sign.
Again we begin with~\eqref{dfXi}, recall~\( \mathfrak{B} = \mleft\{ p_{1}, \ldots, p_{m} \mright\} \).
At primes the difference in the series on the right hand side is
\begin{equation}
  \label{primedf}
  \zeta\mleft( s \mright) - \zeta_{\mathfrak{B}}\mleft( s \mright)
  = \sum_{P^{-}\mleft( n \mright) > p_{m}} \frac{1}{n^{s}}
  = \frac{1}{p_{m + 1}^{s}} + \ldots,
\end{equation}
while at composite numbers~\( k \) it becomes
\begin{align*}
  \zeta\mleft( s \mright) - \zeta_{\mathfrak{B} \cup \mleft\{ k \mright\}}
  & = \zeta\mleft( s \mright) - \mleft( 1 + \frac{1}{2^{s}} + \ldots
  + \frac{1}{\mleft( k - 1 \mright)^{s}} + \frac{2}{k^{s}} + \ldots \mright) \\
  & = - \frac{1}{k^{s}} + \ldots.
  \numberthis \label{compositedf}
\end{align*}
By the inequality~\eqref{pQuality}, the size of~\eqref{dfXi} is majorized by
the first term in the series difference.
Then by the opposite sign in \eqref{primedf}~and~\eqref{compositedf},
we conclude that composites indeed give the opposite sign to primes,
and \autoref{thm:sign} is proved.

\section{Proof of~\autoref{thm:rhsc}}

First we show the backward direction,
the proof of which is similar to~\cite[Theorem~2(i)]{nastasescuClassApproximationsRiemann2022}.
Assume Riemann Hypothesis and Simplicity Conjecture,
and suppose bounded Riemann Hypothesis fails. That is,
there exists a~\(  {l}_0 \in \mathbb{N} \) such that
for all~\( u \), the \(  {l}_0 \)-th~zero of~\( \zeta_{u}^{\approx}\mleft( s \mright) \)
is off the critical line. Let \( \rho_{0} \) be the \(  {l}_0 \)-th
simple zero of the Riemann zeta function.
Let \( \mathcal{C}_{2} \) be a contour that contains
the critical strip up to height~\( \Im \rho_{0} + \varepsilon \),
but omits the next zero.
It also starts from height~\( 2 \) so that the poles~\( s = 0, 1 \) are outside.
Next choose another contour~\( \mathcal{C}_{1} \) around the zero~\( \rho_{0} \)
inside~\( \mathcal{C}_{2} \) and omits the critical line.
Let \( \delta_{1} = \inf \mleft\{ \abs{ \xi\mleft( s \mright) } \colon s \in \mathcal{C}_{1} \mright\} \)
~and~\( \delta_{2} = \inf \mleft\{ \abs{ \xi\mleft( s \mright) } \colon s \in \mathcal{C}_{2} \mright\} \).
By~\cite[Theorem~1]{nastasescuClassApproximationsRiemann2022}, there exists
a~\( u_{0} \) such that for all~\( u > u_{0} \), we have
\( \sup \mleft\{ \abs{ \xi\mleft( s \mright) - \xi_{u}^{\approx}\mleft( s \mright) }
\colon s \in \mathcal{C}_{1} \mright\} < \delta_{1} \)
and \( \sup \mleft\{ \abs{ \xi\mleft( s \mright) - \xi_{u}^{\approx}\mleft( s \mright) }
\colon s \in \mathcal{C}_{2} \mright\} < \delta_{2} \).
Now Rouch\'{e}'s theorem tells us that inside \( \mathcal{C}_{1} \)
~and~\( \mathcal{C}_{2} \), \( \xi\mleft( s \mright) \)~and~
\( \xi_{u}^{\approx}\mleft( s \mright) \) will have the same number
of zeros. Thus \( \xi\mleft( s \mright) \) will have a zero inside~\( \mathcal{C}_{1} \),
contradicting the Riemann Hypothesis.

Now for the forward direction. We only need to show the Bounded Riemann Hypothesis implies
the Simplicity Conjecture, since in~\cite[Theorem~2(i)]{nastasescuClassApproximationsRiemann2022}
it is already proven that Riemann Hypothesis is implied.
We start by assuming the Bounded Riemann Hypothesis and Riemann Hypothesis,
and suppose the Simplicity Conjecture fails. Let \( \rho_{0} = 1/2 + i \gamma_{0} \)
be a zero of the Riemann zeta function of multiplicity~\( k_{0} \ge 2 \).
On the one hand, the series expansion around~\( \rho_{0} \) gives us
\begin{equation}
  \label{taylor}
  \abs{ \xi\mleft( s \mright) - c_{0} \mleft( s - \rho_{0} \mright)^{k_{0}} }
  \le C \abs{ s - \rho_{0} }^{k_{0} + 1},
\end{equation}
for some~\( c_{0} \ne 0 \).
On the other hand, by~\cite[Theorem~1]{nastasescuClassApproximationsRiemann2022} we know
\begin{equation}
  \label{diffAsymp}
  \xi\mleft( s \mright) - \xi_{u}^{\approx}\mleft( s \mright)
  = \frac{2 \mleft( -1 \mright)^{m}
  \mleft( 2 \pi p_{m + 1}^{2} \mright)^{2m + 2} e^{- \pi p_{m + 1}^{2}}}
  {s^{m} \mleft( 1 - s \mright)^{m}}
  a\mleft( 2 \pi p_{m + 1}^{2}, s \mright)
  \mleft( 1 + O\mleft( \frac{1}{\log u} \mright) \mright),
\end{equation}
where
\begin{equation*}
  a\mleft( s_{0}, s \mright)
  = \frac{1}{s_{0} - s} + \frac{1}{s_{0} - 1 + s}
  = \frac{2 s_{0} - 1}{\mleft( s_{0} - s \mright) \mleft( s_{0} - 1 + s \mright)}.
\end{equation*}

Define a new function
\begin{equation*}
  G\mleft( s \mright) \coloneqq
  c_{0} \mleft( s - \rho_{0} \mright)^{k_{0}}
  - \frac{2 \mleft( -1 \mright)^{m}
  \mleft( 2 \pi p_{m + 1}^{2} \mright)^{2m + 2} e^{- \pi p_{m + 1}^{2}}}
  {s^{m} \mleft( 1 - s \mright)^{m}}
  a\mleft( 2 \pi p_{m + 1}^{2}, s \mright).
\end{equation*}
The equation~\( G\mleft( s \mright) = 0 \) is equivalent to
\begin{multline}
  \label{gVanish}
  \mleft( 1 + \frac{z}{\rho_{0}} \mright)^{m}
  \mleft( 1 - \frac{z}{\overline{\rho_{0}}} \mright)^{m}
  \mleft( 1 - \frac{\rho_{0} + z}{2 \pi p_{m + 1}^{2}} \mright)
  \mleft( 1 - \frac{\overline{\rho_{0}} - z}{2 \pi p_{m + 1}^{2}} \mright)
  z^{k_{0}} \\
  = \frac{2}{c_{0}} \mleft( \frac{-1}{1/4 + \gamma_{0}^{2}} \mright)^{m}
  \mleft( 2 \pi p_{m + 1}^{2} \mright)^{2 m}
  e^{- \pi p_{m + 1}^{2}} \mleft( 4 \pi p_{m + 1}^{2} - 1 \mright),
\end{multline}
where we write \( s = \rho_{0} + z \),
and want to show \( G\mleft( s \mright) \) must have \( k_{0} \)~zeros
inside a small disc around~\( \rho_{0} \).
The right hand side of~\eqref{gVanish} is independent of~\( z \),
and it gives us the radius of this disc as
\begin{equation}
  \label{radius}
  r_{m} =
  \mleft( \frac{2}{\abs{ c_{0} }} \mleft( 1/4 + \gamma_{0}^{2} \mright)^{-m}
  \mleft( 2 \pi p_{m + 1}^{2} \mright)^{2 m}
  e^{- \pi p_{m + 1}^{2}} \mleft( 4 \pi p_{m + 1}^{2} - 1 \mright) \mright)^{1 / k_{0}}
  \ll e^{- \frac{\pi m^{2} \log m}{k_{0}} }.
\end{equation}
As~\( m \to \infty \), we know \( p_{m} \sim m \log m \).
In~\eqref{gVanish}, on the one hand we have
\begin{equation*}
  \abs{ \mleft( 1 - \frac{\rho_{0} + z}{2 \pi p_{m + 1}^{2}} \mright)
  \mleft( 1 - \frac{\overline{\rho_{0}} - z}{2 \pi p_{m + 1}^{2}} \mright) }
  \sim 1.
\end{equation*}
On the other, when \( \abs{ z } \le r_{m} \) in~\eqref{radius},
by the elementary inequality
\begin{equation*}
  1 + m \abs{ \frac{z}{\rho_{0}} }
  \le \abs{ 1 + \frac{z}{\rho_{0}} }^{m}
  \le \frac{1}{1 - m \abs{ z / \rho_{0} }},
\end{equation*}
we know
\begin{equation*}
  \abs{ \mleft( 1 + \frac{z}{\rho_{0}} \mright)^{m}
  \mleft( 1 - \frac{z}{\overline{\rho_{0}}} \mright)^{m} }
  \sim 1
\end{equation*}
as well. Thus by~\eqref{gVanish} we see \( G\mleft( s \mright) \) will have
\( k_{0} \)~zeros located at the vertices of a regular \( k_{0} \)-gon
of radius~\( r_{m} \) in~\eqref{radius}
centered at~\( \rho_{0} \), as~\( m \to \infty \).
Denote those zeros as~\( w_{j} \) for~\( j \in \mleft[ 1, k_{0} \mright] \).

Next we want to show that when~\( k_{0} > 2 \),
this contradicts the Bounded Riemann Hypothesis
by implying \( \xi_{u}^{\approx}\mleft( s \mright) \)~having zeros off the critical line.
Note in this case some~\( w_{j} \)  {should}
lie off the critical line.
For a small~\( \delta > 0 \) consider the circle~\( C_{j} \)
centered at an arbitrary zero~\( w_{j} \) of radius~\( \delta r_{m} \).
For~\( s \in C_{j} \) we use the polynomial factorization of~\( G\mleft( s \mright) \)
as~\( m \to \infty \) to show the lower bound
\begin{align*}
  \abs{ G\mleft( s \mright) }
  & \sim \abs{ c_{0} } \prod_{k = 1}^{k_{0}} \abs{ s - w_{k} } \\
  & = \abs{ c_{0} } \delta r_{m} \prod_{k \ne j} \abs{ s - w_{k} } \\
  & \ge \abs{ c_{0} } \delta r_{m}^{k_{0}}
  \prod_{k \ne j} \mleft( \abs{ e^{2 \pi i j / k_{0}} - e^{2 \pi i k / k_{0}} } - \delta \mright),
\end{align*}
since each term in the product is
\begin{align*}
  \abs{ s - w_{k} }
  & \ge \abs{ w_{j} - w_{k} } - \abs{ w_{j} - s } \\
  & = r_{m} \mleft( \abs{ e^{2 \pi i j / k_{0}} - e^{2 \pi i k / k_{0}} } - \delta \mright).
\end{align*}

An application of triangle inequality tells us that
\begin{equation*}
  \abs{ \xi_{u}^{\approx}\mleft( s \mright) - G\mleft( s \mright) }
  \le C \abs{ s - \rho_{0} }^{k_{0} + 1}
  + D \frac{2 \mleft( -1 \mright)^{m}
  \mleft( 2 \pi p_{m + 1}^{2} \mright)^{2m + 2} e^{- \pi p_{m + 1}^{2}}}
  {s^{m} \mleft( 1 - s \mright)^{m} \log u}
  a\mleft( 2 \pi p_{m + 1}^{2}, s \mright),
\end{equation*}
where \( C \)~and~\( D \) are implied constants
in \eqref{taylor}~and~\eqref{diffAsymp} respectively.
The first term is bounded by
\begin{align*}
  C \abs{ s - \rho_{0} }^{k_{0} + 1}
  & \le C \mleft( \abs{ s - w_{j} } + \abs{ w_{j} - \rho_{0} } \mright)^{k_{0} + 1} \\
  & = C \mleft( \mleft( 1 + \delta \mright) r_{m} \mright)^{k_{0} + 1}.
\end{align*}
While by \eqref{gVanish}~and~\eqref{radius},
the second is bounded by
\begin{equation*}
  D \abs{ c_{0} } \frac{r_{m}^{k_{0}}}{\log m},
\end{equation*}
as~\( m \to \infty \).
So for~\( m \) large enough and~\( s \in C_{j} \), we have
\begin{equation*}
  \abs{ \xi_{u}^{\approx}\mleft( s \mright) - G\mleft( s \mright) } \le \abs{ G\mleft( s \mright) }.
\end{equation*}
Then Rouch\'{e}'s theorem implies \( \xi_{u}^{\approx}\mleft( s \mright) \)
will have a zero inside each~\( C_{j} \), which is off the critical line,
contradicting the Bounded Riemann Hypothesis.

When \( k_{0} = 2 \), it is possible that both~\( w_{1}, w_{2} \)
lie on the critical line. We consider the simple zero~\( \rho_{1} \)
closest to but higher than~\( \rho_{0} \).
Such zero pairs~\( \rho_{0}, \rho_{1} \) always exist
due to~\cite{levinsonMoreOneThird1974},
where it is shown that a positive proportion of the zeros are simple.
For large enough~\( u \),
a contour whose inside containing only~\( \rho_{0}, \rho_{1} \)
will have \( 3 \)~zeros of~\( \xi_{u}^{\approx}\mleft( s \mright) \) inside,
counting multiplicity.
By the sign alternation of the zeros in~
\cite[Theorem~2(ii)]{nastasescuClassApproximationsRiemann2022},
we know at least \( 2 \)~of the \( 3 \)~zeros of~\( \xi_{u}^{\approx}\mleft( s \mright) \)
will have to lie off the critical line.
Then a similar application of Rouch\'{e}'s theorem
will contradict the Bounded Riemann Hypothesis in this case as well.

\printbibliography

\end{document}